\documentclass{amsart}

\newcommand{\bba}{\mathbb{A}}\newcommand{\bbf}{\mathbb{F}}

\newcommand{\bbq}{\mathbb{Q}}
\newcommand{\bbp}{\mathbb{P}}\newcommand{\bbz}{\mathbb{Z}}
\newcommand{\caa}{\mathcal{A}}\newcommand{\cab}{\mathcal{B}}

\newcommand{\caf}{\mathcal{F}}
\newcommand{\cak}{\mathcal{K}}
\newcommand{\call}{\mathcal{L}}
\newcommand{\can}{\mathcal{N}}
\newcommand{\cao}{\mathcal{O}}
\newcommand{\cax}{\mathcal{X}}\newcommand{\cay}{\mathcal{Y}}\newcommand{\caz}{\mathcal{Z}}
\newcommand{\frA}{\mathfrak{A}}
\newcommand{\frF}{\mathfrak{F}}

\newcommand{\frO}{\mathfrak{O}}
\newcommand{\frp}{\mathfrak{p}}
\DeclareMathOperator{\aut}{Aut}

 \DeclareMathOperator{\et}{\text{\'et}}\DeclareMathOperator{\Frob}{Frob}
\DeclareMathOperator{\gal}{Gal}\DeclareMathOperator{\gl}{GL}
\DeclareMathOperator{\ns}{NS}

\newcommand{\ov}{\overline}\DeclareMathOperator{\pic}{Pic}
\DeclareMathOperator*{\res}{Res}
\DeclareMathOperator{\rk}{rank}\DeclareMathOperator{\spec}{Spec}
\DeclareMathOperator{\tr}{Tr}
\numberwithin{equation}{section}
    \newtheorem{theorem}{Theorem}
\numberwithin{theorem}{section} \theoremstyle{plain}

\newtheorem{conjecture}[theorem]{Conjecture}

\newtheorem{lemma}[theorem]{Lemma}

\newtheorem{proposition}[theorem]{Proposition}

\theoremstyle{definition}
\newtheorem{definition}[theorem]{Definition}

\newtheorem{remark}[theorem]{Remark}

\begin{document}
\title[Rank of Jacobian varieties]{On the variation of the rank of Jacobian varieties on unramified abelian towers\\ over
number fields}

\author{Am\'{\i}lcar Pacheco}

\address{Universidade Federal do Rio de Janeiro (Universidade do Brasil)\\
Departamento de Ma\-te\-m\'a\-ti\-ca Pura\\
Rua Guai\-aquil 83, Cachambi, 20785-050 Rio de Janeiro, RJ, Brasil}

\email{amilcar@impa.br}

\thanks{This work was partially supported by CNPq research grant 300896/91-3, Pronex 41.96.0830.00 and CNPq Edital
Universal 470099/2003-8.}

\date{November 15, 2003}

\begin{abstract}
Let $C$ be a smooth projective curve defined over a number field $k$, $X/k(C)$ a smooth projective curve of positive
genus, $J_X$ the Jacobian variety of $X$ and $(\tau,B)$ the $k(C)/k$-trace of $J_X$. We estimate how the rank of
$J_X(k(C))/\tau B(k)$ varies when we take an unramified abelian cover $\pi:C'\to C$ defined over $k$.
\end{abstract}

\maketitle

\section{Introduction}

Let $C$ be a smooth projective curve defined over a number field $k$,
$K=k(C)$ its function field and $A$ an abelian variety defined over $K$. Let $(\tau,B)$ be the $K/k$-trace of $A$ and $\ov{k}$ an algebraic closure of $k$. A theorem of
N\'eron and Lang states that both groups $A(K)/\tau B(k)$ and $A(\ov{k}(C))/\tau B(\ov{k})$ are finitely generated. In
this paper we will consider the rank of these two groups in the case where $A$ is the Jacobian variety $J_X$ of a smooth projective curve $X$
defined over $K$ of genus $g_X\ge1$.

Let $g_C$ be the genus of $C$, $\frF_{J_X}$ the conductor divisor of $J_X$ and $f_{J_X}$ its degree. Ogg proves in
\cite[VI, p. 19]{ogg67} the following geometric upper bound
\begin{equation}\label{ogg1}
\rk\left(\frac{J_X(\ov{k}(C))}{\tau B(\ov{k})}\right)\le2g_X(2g_C-2)+f_{J_X}+4\dim(B).
\end{equation}
(See also \cite[Theorem 2]{ogg62}). In the case where $X=E$ is an elliptic curve, (\ref{ogg1}) reduces to
\begin{equation}\label{shi1}
\rk(E(\ov{k}(C)))\le4g_C-4+f_E,
\end{equation}
which is a result due to Shioda \cite[Corollary 2]{shi92}.

The goal of this paper is to study the variation of the rank of $J_X(K)/\tau B(k)$ under an unramified abelian covering
$\pi:C'\to C$ defined over $k$ generalizing a result of Silverman
\cite[Theorem 13]{sil02} proved in the case where $X$ was an elliptic
curve.

Let $\cax$ be a model of $X$, i.e., $\cax$ is a smooth projective surface defined over $k$ and $X$ is the generic fiber
of the proper flat morphism $\phi:\cax\to C$ also defined over $k$ of relative dimension 1. Let $K'=k(C')$,
$\cax'=\cax\times_CC'$, $\phi':\cax'\to C'$ the morphism obtained from $\phi$ by extending the base $C$ to $C'$ via
$\pi$. Its generic fiber $X'$ is isomorphic to $X\times_KK'$, so the Jacobian variety $J_{X'}$ of $X'$ is isomorphic to
$J_X\times_KK'$. As a consequence the genus $g_{X'}$ of $X'$ equals $g_X$. Observe that since the extension $K'/K$ is
geometric the $K/k$-trace $(\tau,B)$ of $J_X$ is $k$-isomorphic to the $K'/k$-trace of $J_{X'}$. Let $\frF_{J_{X'}}$ be
the conductor divisor of $J_{X'}$ and $f_{J_{X'}}=\deg(\frF_{J_{X'}})$. The geometric bound (\ref{ogg1}) applied to
$C'$ gives
\begin{equation}\label{arbd}
\rk\left(\frac{J_{X'}(K')}{\tau B(k)}\right)\le\rk\left(\frac{J_{X'}(\ov{k}(C'))}{\tau
B(\ov{k})}\right)\le2g_X(2g_{C'}-2)+f_{J_{X'}}+4\dim(B).
\end{equation}

In order to state our main theorem we recall Tate's conjecture for smooth projective surfaces $\cax$ defined over a
number field $k$. Let $l$ be a perfect field and $\ov{l}$ an algebraic closure of $l$. For each algebraic variety
$\cay$ defined over $l$ we denote by $H^i(\cay)$ its $i$-th \'etale cohomology group
$H^i_{\et}(\cay\times_l\ov{l},\bbq_{\ell})$.

Let $\cao_k$ be the ring of integers of $k$ and $\frp$ a prime ideal of $\cao_k$. Let $\Frob_{\frp}\in G_k$ be a
Frobenius automorphism corresponding to $\frp$ and $I_{\frp}\subset G_k$ the inertia group of $\frp$ (well defined up
to conjugation). Define the $L$-function
$$
L_2(\cax/k,s)=\prod_{\frp}\det(1-\Frob_{\frp}q_{\frp}^{-s}\,|\,H^2(\cax)^{I_{\frp}})^{-1}.
$$
Let $\pic(\cax)$ be the Picard group of $\cax$, $\pic^0(\cax)$ the subgroup of divisors algebraically equivalent to
zero, $\ns(\cax)=\pic(\cax)/\pic^0(\cax)$ the N\'eron-Severi group of $\cax$ and $\ns(\cax/k)$ the subgroup of divisor
classes of $\ns(\cax)$ which are defined over $k$. The group $\ns(\cax)$ is finitely generated, hence the same holds
for $\ns(\cax/k)$.

\begin{conjecture}(Tate's conjecture, \cite[Conjecture 2]{tate65})\label{tateconj}
$L_2(\cax/k,s)$ has a pole at $s=2$ of order $\rk(\ns(\cax/k))$.
\end{conjecture}

\begin{remark}
\begin{enumerate}
\item This is a special case of Tate's conjecture which concerns algebraic varieties and algebraic cycles.

\item We do not need to use the hypothesis of the existence of a meromorphic continuation of $L_2(\cax/k,s)$ to
$\Re(s)=2$ interpreting the sentence ``$L_2(\cax/k,s)$ has a pole of order $t$ in $s=2$'' meaning
$$
\lim_{\Re(s)>2,s\to2}(s-2)^tL_2(\cax/k,s)=\beta\ne0.
$$
Moreover, if $f(s)$ is a holomorphic function in $\Re(s)>\lambda$ and $\lim_{s\to\lambda}(s-\lambda)f(s)=\beta\ne0$, we
will call $\beta$ the residue of the function $f(s)$ in $s=\lambda$ and we will write $\res_{s=\lambda}(f(s))=\alpha$.
\end{enumerate}
\end{remark}

Let $G_k=\gal(\ov{k}/k)$, $\caa=\aut_{\ov{k}}(C'/C)$ the subgroup of the group $\aut_{\ov{k}}(C')$ of
$\ov{k}$-automorphisms of $C'$ which fixes the points of $C$ and $\frO_{G_k}(\caa)$ the set of $G_k$-orbits of $\caa$.

\begin{theorem}\label{mainthm}
Suppose Tate's conjecture is true for $\cax'/k$. Then
\begin{equation}\label{rk1}
\rk\left(\frac{J_{X'}(K')}{\tau B(k)}\right)\le\frac{\#\frO_{G_k}(\caa)}{|\caa|}(2g_X(2g_{C'}-2)+f_{J_{X'}}).
\end{equation}
\end{theorem}

In the case where $X=E$ is an elliptic curve, (\ref{rk1}) is indeed refined to \cite[Theorem 1]{sil02}
\begin{equation}\label{rkec1}
\rk(E'(K'))\le\frac{\#\frO_{G_k}(\caa)}{|\caa|}(4g_{C'}-4+f_{E'}).
\end{equation}
In particular, (\ref{rk1}) (resp. (\ref{rkec1})) will only give an improvement of (\ref{ogg1}) (resp. (\ref{shi1})) if
the action of $G_k$ on $\caa$ is non trivial.

This improvement actually happens in at least two instances. First when $C=C'$ is an elliptic curve and $\pi$ is equal to the
multiplication by an integer $n\ge1$ map in $C$, then $\caa$ is the subgroup of $n$-torsion points $C[n]$ of $C$. A
theorem of Serre \cite{serre72} states that this action is highly non-trivial as $n$ grows. We also consider the case
in which $C'$ is the pullback of $C$ under the multiplication by $n$ map in the Jacobian variety $J_C$ of $C$ and $\pi$
is the corresponding unramified abelian covering. In this case $\caa$ is the subgroup of $n$-torsion points $J_C[n]$ of
$J_C$. Under the hypothesis that $k$ is sufficiently large, Serre
\cite{serre86rib} extended the previous result from elliptic curves to
abelian varieties. As a consequence, the action of $G_k$ on $J_C[n]$ is also highly non-trivial as $n$ grows.

In \S\,2 we describe the connection between Tate's conjecture and the generalized analytic Nagao's conjecture. In \S\,3
we use Deligne's equidistribution theorem to give an upper bound for
the absolute value $|\frA_{\frp}(\cax)|$ of the average trace of
Frobenius in terms of the degree of
the conductor $\frF_{J_X}$. We also show that the conductor behaves well with respect to finite unramified base
extensions. In \S\,4 we obtain auxiliary results through
counting rational points. In \S\,5 we prove Theorem \ref{mainthm}. In \S\,6 we give applications analyzing the variation
of the rank in special unramified abelian towers over number fields. In these special cases, we show that the rank will
grow more slowly along the tower than the geometric bound (\ref{rk1}).

\section{Tate's conjecture and the generalized Nagao's conjecture}

Given a prime ideal $\frp$ of $\cao_k$ and an algebraic variety $\cay$ defined over $k$, we will denote by
$\tilde{\cay}_{\frp}$ its reduction modulo $\frp$. Given an algebraic variety $\cay$ defined over a perfect field $l$,
let $H^i_{c}(\cay)$ be its $i$-th cohomology group $H^i_{c}(\cay\times_l\ov{l},\bbq_{\ell})$ with compact support.

Let $S$ be a finite set of prime ideals of $\cao_k$ (which will be enlarged as needed). First we assume that for every
$\frp\notin S$, $\tilde{\cax}_{\frp}$ (resp. $\tilde{C}_{\frp}$) is a smooth projective surface (resp. curve) over the
residue field $\bbf_{\frp}$ of $\frp$ of cardinality $q_{\frp}$ and that the reduction
$\tilde{\phi}_{\frp}:\tilde{\cax}_{\frp}\to\tilde{C}_{\frp}$ of $\phi$ modulo $\frp$ is a proper flat morphism of
relative dimension 1 defined over $\bbf_{\frp}$.

For each $z\in\tilde{C}_{\frp}(\bbf_{\frp})$, let $\tilde{\cax}_{\frp,z}=\tilde{\phi}_{\frp}^{-1}(z)$ be the fiber of
$\tilde{\phi}_{\frp}$ at $z$. Let ${F}_{\frp}$ be the topological generator of $\gal(\ov{\bbf}_{\frp}/\bbf_{\frp})$.
Denote also by $F_{\frp}$ its induced automorphism on $H^1(\tilde{\cax}_{\frp,z})$ (resp.
$H^1_{c}(\tilde{\cax}_{\frp,z})$). Let $\Delta=\{z\in C\,|\,\phi^{-1}(z)$ is not smooth$\}$ be the discriminant locus
of $\phi$. After discarding a finite number of prime ideals $\frp$ of $\cao_k$, we may assume that for every
$\frp\notin S$ the discriminant locus $\tilde{\Delta}_{\frp}$ of $\tilde{\phi}_{\frp}$ is equal to reduction modulo
$\frp$ of $\Delta$.

For every $z\in(\tilde{C}_{\frp}-\tilde{\Delta}_{\frp})(\bbf_{\frp})$, let
$a_{\frp}(\tilde{\cax}_{\frp,z})=\tr({F}_{\frp}\,|\,H^1(\tilde{\cax}_{\frp,z}))$ and for every
$z\in\tilde{\Delta}_{\frp}(\bbf_{\frp})$ let
$a_{\frp}(\tilde{\cax}_{\frp,z})=\tr(F_{\frp}\,|\,H^1_{c}(\tilde{\cax}_{\frp,z})$. The \emph{average trace of
Frobenius} is defined by
$$
\frA_{\frp}(\cax)=\frac1{q_{\frp}}\sum_{z\in\tilde{C}_{\frp}(\bbf_{\frp})}a_{\frp}(\tilde{\cax}_{\frp,z}).
$$

 Let $a_{\frp}(B)=\tr(\Frob_{\frp}\,|\,H^1(B)^{I_{\frp}})$. By base change (cf. \cite{mil80} or \cite[Appendix C]{har86})
 this number equals
 $\tr(F_{\frp}\,|\,H^1(\tilde{B}_{\frp}))$. The
\emph{reduced average trace of Frobenius} is defined by
$$
\frA_{\frp}^*(\cax)=\frA_{\frp}(\cax)-a_{\frp}(B).
$$

\begin{theorem}\label{thmhp}\cite[Th\'eor\`eme 1.3]{hinpac03}
Tate's Conjecture \ref{tateconj} implies the generalized analytic Nagao's conjecture:
$$
\res_{s=1}\left(\sum_{\frp\notin S}-\frA_{\frp}^*(\cax)\frac{\log
q_{\frp}}{q_{\frp}^s}\right)=\rk\left(\frac{J_X(K)}{\tau B(k)}\right).
$$
\end{theorem}

\section{Equidistribution theorem and the conductor}

\begin{remark}\label{remfrob}
It follows from Weil's theorem (the Riemann hypothesis for curves over finite fields) that for every
$z\in(\tilde{C}_{\frp}-\tilde{\Delta}_{\frp})(\bbf_{\frp})$ all eigenvalues of $F_{\frp}$ acting on
$H^1(\tilde{\cax}_{\frp,z})$ have absolute value $q_{\frp}^{1/2}$, thus
$a_{\frp}(\tilde{\cax}_{\frp,z})=O(q_{\frp}^{1/2})$. For $z\in\tilde{\Delta}_{\frp}(\bbf_{\frp})$ it is a result due to
Deligne \cite[Th\'eor\`eme 3.3.1]{del81} that all the eigenvalues of $F_{\frp}$ acting on
$H^1_{c}(\tilde{\cax}_{\frp,z})$ have absolute value at most $q_{\frp}^{1/2}$, so once again
$a_{\frp}(\tilde{\cax}_{\frp,z})=O(q_{\frp}^{1/2})$. As a consequence,
$\frA_{\frp}(\cax)=O(q_{\frp}^{1/2})$. Theorem \ref{avgtr1} uses Deligne's equidistribution theorem and local monodromy to improve this estimate (cf. \cite[Th\'eor\`eme
3.5.3]{del81} or \cite[(3.6.3)]{kat88}). A similar improvement in the case where $X$ is an elliptic curve was obtained
by Silverman \cite[Theorem 6]{sil02}.
\end{remark}

Let $\caz$ be a smooth projective curve defined over a perfect
field $l$ of characteristic $p\ge0$, $\call=l(\caz)$ its function
field and $A$ be an abelian variety $\call$. Let $\ell\ne p$ be a prime number, $A[\ell]$ the subgroup of $\ell$-torsion
points of $A$, $T_{\ell}(A)$ the $\ell$-adic Tate module of $A$ and
$V_{\ell}(A)=T_{\ell}(A)\otimes_{\bbz_{\ell}}\bbq_{\ell}$. Let
$\call^s$ be a separable closure of $\call$. For every place $v$ of
$\call$ denote by $I_v\subset\gal(\call^s/\call)$ the inertia subgroup corresponding to $v$ (which is well defined
up to conjugation).

\begin{definition}\label{defcond}
The multiplicity of the conductor $\frF_A$ of $A$ at $v$ is equal to a sum of two numbers, the tame part $\epsilon_v$
of $\frF_A$ at $v$ and the wild part $\delta_v$ of $\frF_A$ in $v$. The first number is defined as
$\epsilon_v=\text{codim}(V_{\ell}(A)^{I_v})$, where $V_{\ell}(A)^{I_v}$ denotes the set of elements fixed by the action
of $I_v$. If $N(A)_v$ denotes the N\'eron model of $A$ over $\spec(\cao_{v})$
and $\can(A)_v^0$ is the connected component of the special fiber
$\can(A)_v$ of $N(A)_v$, let $u_v$ (resp. $r_v$) be the unipotent (resp.
reducible) rank of $\can(A)^0_v$, then $\epsilon_v=u_v+2r_v$.

The second number is defined as follows. Let $k_v$ be the residue field of $v$ and $l_v/k_v$ a finite
Galois extension such that $G_{l_v}=\gal(\ov{l}_v/l_v)$ acts trivially on $A[\ell]$. Let
$G_v=\gal(l_v/k_v)$, so $A[\ell]$ can be regarded as a $G_v$-module. Let $P_v$ be a projective
$\bbz_{\ell}[G_v]$-module whose character is the Swan character of $G_v$. Then
$\delta_v:=\dim_{\bbf_{\ell}}(\text{Hom}_{\bbz_{\ell}[G_v]}(P_v,A[\ell]))$
(cf. \cite[\S1]{ogg67}). This is a non-negative 
integer and is in fact independent from the choice of $l_v$.
\end{definition}

Suppose $\frp\notin S$. Let $\tilde{U}_{\frp}:=\tilde{C}_{\frp}-\tilde{\Delta}_{\frp}$,
$\caf_{\frp}=R^1(\tilde{f}_{\frp})_*\bbq_{\ell}$ and 
$H^1_{c}(\tilde{U}_{\frp}\times_{\bbf_{\frp}}\ov{\bbf}_{\frp},\caf_{\frp})^{\text{wt}\le0}$ the part of
weight $\le0$ of $H^1_{c}(\tilde{U}_{\frp}\times_{\bbf_{\frp}}\ov{\bbf}_{\frp},\caf_{\frp})$ (which is mixed of
weight $\le1$). Let $\ov{\eta}$ be the geometric generic point of
$\tilde{U}_{\frp}\times_{\bbf_{\frp}}\ov{\bbf}_{\frp}$ and
$\tilde{X}_{\frp}$ the generic fiber of $\tilde{f}_{\frp}$. Then by proper
base change $(\caf_{\frp})_{\ov{\eta}}\cong H^1(\tilde{X}_{\frp})\cong
H^1(\pic^0(\tilde{X}_{\frp}))\cong V_{\ell}(\pic^0(\tilde{X}_{\frp}))^{\vee}$.

\begin{lemma}\label{wtcoh}
$$
H^1_{c}(\tilde{U}_{\frp}\times_{\bbf_{\frp}}\ov{\bbf}_{\frp},\caf_{\frp})^{\text{wt}\le0}\cong\bigoplus_{z\in\tilde{\Delta}_{\frp}(\bbf_{\frp})}((\caf_{\frp})_{\ov{\eta}})^{I_z}\text{
  and }
$$
$$
\text{codim}(H^1_{c}(\tilde{U}_{\frp}\times_{\bbf_{\frp}}\ov{\bbf}_{\frp},\caf_{\frp})^{\text{wt}\le0})\le
2g_X(2g_C-2)+\sum_{z\in\tilde{\Delta}_{\frp}(\bbf_{\frp})}\epsilon_z.
$$
\end{lemma}

\begin{proof}
The first statement follows from \cite[Lemme 4.1]{mic95} replacing
$\caf$ by $\caf_{\frp}$, $U_p$ by
$\tilde{U}_{\frp}$ and $\bbp^1$ by $\tilde{C}_{\frp}$. For the second
statement it follows from the first isomorphism that 
$$
\begin{aligned}
&\dim(H^1_{c}(\tilde{U}_{\frp}\times_{\bbf_{\frp}}\ov{\bbf}_{\frp},\caf)^{\text{wt}\le0})=\sum_{z\in\tilde{\Delta}_{\frp}(\bbf_{\frp})}\dim((V_{\ell}(\pic^0(\tilde{X}_{\frp}))^{\vee})^{I_c})\\
&=\sum_{z\in\tilde{\Delta}_{\frp}(\bbf_{\frp})}(2g_{\tilde{X}_{\frp}}-\epsilon_c).
\end{aligned}
$$
Since
$$
\dim(H^1_{c}(\tilde{U}_{\frp}\times_{\bbf_{\frp}}\ov{\bbf}_{\frp},\caf_{\frp}))=2g_{\tilde{X}_{\frp}}(2g_{\tilde{C}_{\frp}}-2+s_{J_{\tilde{X}_{\frp}}})=2g_{\tilde{X}_{\frp}}(2g_{\tilde{C}_{\frp}}-2)+\sum_{z\in\tilde{\Delta}_{\frp}(\bbf_{\frp})}(2g_{\tilde{X}_{\frp}}),
$$
where $s_{J_{\tilde{X}_{\frp}}}$ denotes
$\#\text{supp}(\frF_{J_{\tilde{X}_{\frp}}})$, we conclude that
\begin{align*}
&\text{codim}(H^1_{c}(\tilde{U}_{\frp}\times_{\bbf_{\frp}}\ov{\bbf}_{\frp},\caf)^{\text{wt}\le0})\\
&=2g_{\tilde{X}_{\frp}}(2g_{\tilde{C}_{\frp}}-2)+\sum_{z\in\tilde{\Delta}_{\frp}(\bbf_{\frp})}(2g_{\tilde{X}_{\frp}})-\sum_{z\in\tilde{\Delta}_{\frp}(\bbf_{\frp})}(2g_{\tilde{X}_{\frp}}-\epsilon_c)\\
&=2g_{\tilde{X}_{\frp}}(2g_{\tilde{C}_{\frp}}-2)+\sum_{z\in\tilde{\Delta}_{\frp}(\bbf_{\frp})}\epsilon_c.
\end{align*}
The lemma now follows from observing that the curve $\tilde{C}_{\frp}$ obtained from $C$ by reducing it modulo $\frp$ has also genus $g_C$. Furthermore,
enlarging the set $S$, if necessary, we may assume that for every $\frp\notin S$ the generic fiber $\tilde{X}_{\frp}$
of $\tilde{\phi}_{\frp}$ equals the reduction  of $X$ modulo $\frp$. In particular, its genus equals $g_X$.
\end{proof}

\begin{theorem}\label{avgtr1}
$$
|\frA_{\frp}(\cax)|\le2g_X(2g_C-2)+f_{J_X}+O(q_{\frp}^{-1/2}).
$$
\end{theorem}

\begin{proof}
By the previous choice of the set $S$ we assume that for every $\frp\notin S$ the generic fiber $\tilde{X}_{\frp}$
of $\tilde{\phi}_{\frp}$ equals the reduction  of $X$ modulo
$\frp$. As a consequence, the
conductor of the Jacobian variety of $\tilde{X}_{\frp}$ equals the reduction of $\frF_{J_X}$ modulo $\frp$. Hence, it
has degree $f_{J_X}$ and $s_{J_X}=\#\text{supp}(\frF_{J_X})=s_{J_{\tilde{X}_{\frp}}}$.

Let
\begin{equation}\label{eqeq1}
\frA'_{\frp}(\cax)=\frac1{\#(\tilde{C}_{\frp}-\tilde{\Delta}_{\frp})(\bbf_{\frp})}
\sum_{z\in(\tilde{C}_{\frp}-\tilde{\Delta}_{\frp})(\bbf_{\frp})}a_{\frp}(\tilde{\cax}_{\frp,z}).
\end{equation}
Deligne's equidistribution theorem (\cite[Th\'eor\`eme 3.5.3]{del81} or \cite[(3.6.3)]{kat88}) states
\begin{equation}\label{eqeq2}
|\frA'_{\frp}(\cax)|\le2g_X(2g_C-2+s_{J_X})+O(q_{\frp}^{-1/2}).
\end{equation}

By Lemma \ref{wtcoh} and the Grothendieck-Lefschetz formula, (\ref{eqeq2}) is refined as
\begin{equation}\label{refcond}\begin{aligned}
|\frA'_{\frp}(\cax)|&\le2g_X(2g_C-2)+\sum_{z\in\tilde{\Delta}_{\frp}(\bbf_{\frp})}\epsilon_z+O(q_{\frp}^{-1/2})\\
&\le2g_X(2g_C-2)+\sum_{v}(\epsilon_v+\delta_v)\deg(v)+O(q_{\frp}^{-1/2})\\
&=2g_X(2g_C-2)+f_{J_{\tilde{X}_{\frp}}}+O(q_{\frp}^{-1/2}),
\end{aligned}\end{equation}
where $v$ runs through the places of $\bbf_{\frp}(\tilde{C}_{\frp})$.
But by our choice of $\frp$ we have $f_{J_{\tilde{X}_{\frp}}}=f_{J_X}$.

By \cite[Th\'eor\`eme 3.3.1]{del81}
\begin{equation}\label{eqeq3}
\frac1{q_{\frp}}\sum_{z\in\tilde{\Delta}_{\frp}(\bbf_{\frp})}a_{\frp}(\tilde{\cax}_{\frp,z})=O(q_{\frp}^{-1/2}).
\end{equation}
Let $a_{\frp}(C)=\tr(\Frob_{\frp}\,|\,H^1(C)^{I_{\frp}})=\tr(F_{\frp}\,|\,H^1(\tilde{C}_{\frp}))$. By Weil's theorem
$$
\#\tilde{C}_{\frp}(\bbf_{\frp})=q_{\frp}+1-a_{\frp}(C).
$$
Observe that
\begin{equation}\label{eqeq4}\begin{aligned}
&\left|\frA'_{\frp}(\cax)-\frac1{q_{\frp}}\sum_{z\in(\tilde{C}_{\frp}-\tilde{\Delta}_{\frp})(\bbf_{\frp})}
a_{\frp}(\tilde{\cax}_{\frp,z})\right|\\
&=\frac{|a_{\frp}(C)-1+\#\tilde{\Delta}_{\frp}(\bbf_{\frp})|}{q_{\frp}
\#(\tilde{C}_{\frp}-\tilde{\Delta}_{\frp})(\bbf_{\frp})}\left|\sum_{z\in(\tilde{C}_{\frp}-\tilde{\Delta}_{\frp})
(\bbf_{\frp})}a_{\frp}(\tilde{\cax}_{\frp,z})\right|\\
&\le\frac{|a_{\frp}(C)-1+\#\tilde{\Delta}_{\frp}(\bbf_{\frp})|}{q_{\frp}}\,O(q_{\frp}^{1/2})=O(1).
\end{aligned}\end{equation}
The theorem now follows from (\ref{refcond}), (\ref{eqeq3}) and (\ref{eqeq4}).
\end{proof}

The conductor behaves well with respect to finite unramified base extensions.

\begin{proposition}(Cf. \cite[Proposition 8]{sil02})\label{condext}
Suppose $\pi:C'\to C$ is unramified (not necessarily abelian). Then
\begin{enumerate}
\item $\frF_{J_{X'}}=\pi^*\frF_{J_X}$.

\item Let $\cak_C$, resp. $\cak_{C'}$, be the canonical divisor of $C$, resp. $C'$, these conductors are related by
$\cak_{C'}=\pi^*\cak_C$.

\item In particular,
\begin{equation}\label{condextf}
2g_{C'}-2+f_{J_{X'}}=|\caa|(2g_C-2+f_{J_X}).
\end{equation}
\end{enumerate}
\end{proposition}

\begin{proof}
Let $v'$ be a place of $K'$ lying over a place $v$ of $K$ via $\pi$. Since $\pi$ is unramified and $k_{v'}/k_v$ is separable (because these fields
are number fields), then by \cite[Chapter 7, Theorem 1, p. 176]{ray90} the morphism
$N(J_X\times_KK_v)\times_{\spec(\cao_{v})}\spec(\cao_{v'})\to N(J_{X'}\times_{K'}K'_{v'})$ of
N\'eron models is an isomorphism. In particular, at the level of the connected components of the special fibers, we
have an isomorphism
$\can(J_X\times_KK_v)^0\times_{\spec(k_v)}\spec(k_{v'})\to\can(J_{X'}\times_{K'}K'_{v'})^0$.
Therefore the tame parts $\epsilon_{v'}$ (resp. $\epsilon_v$) of $\frF_{J_{X'}}$ (resp. $\frF_{J_X}$) at $v'$ (resp.
$v$) are equal. Furthermore, since for every place $v$ of $K$
(resp. $v'$ of $K'$) the residue field $k_v$ (resp.
$k_{v'}$) is a number field, then the extension $K(J_X[\ell])/K$ (resp.
$K'(J_{X'}[\ell])/K'$) is tame. Hence, neither $\frF_{J_X}$ nor $\frF_{J_{X'}}$ have any wild part at
$v$ (resp. $v'$). In particular, $\frF_{J_{X'}}=\pi^*(\frF_{J_X})$. Item (2) follows from \cite[Proposition
IV.2.3]{har86} and (3) follows from (1) and (2) and \cite[IV, 1.3.3]{har86}.
\end{proof}

\section{Counting points}

\begin{lemma}\label{trfrob}
Let $\frp\notin S$ and $\tilde{\pi}_{\frp}:\tilde{C}'_{\frp}\to\tilde{C}_{\frp}$ be the reduction of $\pi$ modulo
$\frp$, $z'\in\tilde{C}'_{\frp}(\bbf_{\frp})$ and $z=\tilde{\pi}_{\frp}(z')\in\tilde{C}_{\frp}(\bbf_{\frp})$. Then
$a_{\frp}(\tilde{\cax}'_{\frp,z'})=a_{\frp}(\tilde{\cax}_{\frp,z})$.
\end{lemma}

\begin{proof}
It suffices to note that the residue field extension $k_{z'}/k_z$ is finite and that
$\tilde{\cax'}_{\frp,z'}\cong\tilde{\cax}_{\frp,z}\times_{\spec(k_z)}\spec(k_{z'})$.
\end{proof}

Let $\cab$ be a subgroup of $\caa$ and $C_{\cab}=C'/\cab$ the intermediate curve $C'\to C_{\cab}\to C$ with
$\aut_{\ov{k}}(C'/C_{\cab})=\cab$. Since $\caa$ is abelian, the cover $C_{\cab}\to C$ is also Galois. Let $L/k$ be a
finite Galois extension sufficiently large so that all intermediate curves $C_{\cab}$ are defined over $L$. If
$\gal(L/k)$ acts on $\cab$, then $C_{\cab}$ is also defined over $k$.

After adding finitely many prime ideals to $S$, we may also assume that for every $\frp\notin S$,
$\tilde{\pi}_{\frp}:\tilde{C}'_{\frp}\to\tilde{C}_{\frp}$ is unramified abelian and
$\tilde{\caa}_{\frp}=\aut_{\ov{\bbf}_{\frp}}(\tilde{C}'_{\frp}/\tilde{C}_{\frp})$ is isomorphic to $\caa$. Let
$\tilde{\bba}_{\frp}=\aut_{\bbf_{\frp}}(\tilde{C}'_{\frp}/\tilde{C}_{\frp})$.

\begin{proposition}\cite[Proposition 11]{sil02}\label{countpt}
Let $\sigma$ be a generator of $\gal(\ov{\bbf}_{\frp}/\bbf_{\frp})$ and $\cab\subset\caa$ the subgroup defined by
$\cab=\{\sigma(a)\circ a^{-1}\,;\,a\in\caa\}$.
\begin{enumerate}
\item The group $\cab$ is defined over $\bbf_{\frp}$, hence the curve $\tilde{C}_{\cab,\frp}=\tilde{C}'_{\frp}/\cab$ is
also defined over $\bbf_{\frp}$.

\item The image of $\tilde{C}'_{\frp}(\bbf_{\frp})$ in $\tilde{C}_{\frp}(\bbf_{\frp})$ coincides with the image of
$\tilde{C}_{\cab,\frp}(\bbf_{\frp})$ in $\tilde{C}_{\frp}(\bbf_{\frp})$.

\item For every $z\in\tilde{C}_{\frp}(\bbf_{\frp})$ such that $z=\tilde{\pi}_{\frp}(z')$ with $z'\in
\tilde{C}'_{\frp}(\bbf_{\frp})$ there exist exactly $|\tilde{\bba}_{\frp}|$ points in $\tilde{C}'_{\frp}(\bbf_{\frp})$
lying over $z$.

\item For every $z\in\tilde{C}_{\frp}(\bbf_{\frp})$ such that $z=\tilde{\pi}_{\cab,\frp}(z_{\cab})$ with
$z_{\cab}\in\tilde{C}_{\cab,\frp}(\bbf_{\frp})$ there exist exactly $(\caa:\cab)$ points of
$\tilde{C}_{\cab,\frp}(\bbf_{\frp})$ lying over $z$.
\end{enumerate}
\end{proposition}

Let $\cab$ be as in Proposition \ref{countpt}. Let $\phi_{\cab}:\cax_{\cab}\to C_{\cab}$ obtained from $\phi$ by the
extension $C_{\cab}\to C$ of the base $C$, $X_{\cab}/k(C_{\cab})$ the generic fiber of $\phi_{\cab}$ and $J_{X_{\cab}}$
its Jacobian variety. We enlarge $S$ even further so that the for every $\frp\notin S$, $\tilde{\cax}_{\cab,\frp}$ and
$\tilde{C}_{\cab,\frp}$ are smooth and $\tilde{\phi}_{\cab,\frp}:\tilde{\cax}_{\cab,\frp}\to\tilde{C}_{\cab,\frp}$ is a proper flat morphism of
relative dimension 1 defined over $\bbf_{\frp}$. Let $\tilde{X}_{\cab,\frp}$ be the generic fiber of
$\tilde{\phi}_{\cab,\frp}$ and $J_{\tilde{X}_{\cab,\frp}}$ its Jacobian variety. Furthermore, we also assume that for
every $\frp\notin S$ the conductor $\frF_{J_{\tilde{X}_{\cab,\frp}}}$ of $J_{\tilde{X}_{\cab,\frp}}$ is equal to the
reduction modulo $\frp$ of the conductor $\frF_{J_{X_{\cab}}}$ of $J_{X_{\cab}}$.

\begin{proposition}(Cf. \cite[Proposition 12]{sil02})\label{avgtr2}
$$
|\frA_{\frp}(\cax')|\le\frac{|\tilde{\bba}_{\frp}|}{|\caa|}(2g_X(2g_{C'}-2)+f_{J_{X'}})+O(q_{\frp}^{-1/2}).
$$
\end{proposition}

\begin{proof}
Let $\cab=\{\sigma(a)\circ a^{-1}\,|\,a\in\caa\}$ be the subgroup defined in Proposition \ref{countpt}. To ease
notations denote $C_{\cab}=C''$ and $\cax_{\cab}=\cax''$. By Lemma \ref{trfrob} applied to
$\tilde{\pi}_{\frp}:\tilde{C}'_{\frp}\to\tilde{C}_{\frp}$ (resp.
$\tilde{\pi}'_{\frp}:\tilde{C}''_{\frp}\to\tilde{C}_{\frp}$) we have
$a_{\frp}(\tilde{\cax}'_{\frp,z'})=a_{\frp}(\tilde{\cax}_{\frp,z})$ (resp.
$a_{\frp}(\tilde{\cax}''_{\frp,z''})=a_{\frp}(\tilde{\cax}_{\frp,z})$) for $z'\in\tilde{C}'_{\frp}(\bbf_{\frp})$ (resp.
$z''\in\tilde{C}''_{\frp}(\bbf_{\frp})$) such that $\tilde{\pi}_{\frp}(z')=z\in\tilde{C}_{\frp}(\bbf_{\frp})$ (resp.
$\tilde{\pi}'_{\frp}(z'')=z\in\tilde{C}_{\frp}(\bbf_{\frp})$). Thus, by Proposition \ref{countpt},
\begin{align*}
q_{\frp}\frA_{\frp}(\cax')&=\sum_{z'\in\tilde{C}'_{\frp}(\bbf_{\frp})}a_{\frp}(\tilde{\cax}'_{\frp,z'})=
|\tilde{\bba}_{\frp}|\sum_{z\in\tilde{\pi}_{\frp}(\tilde{C}'_{\frp}(\bbf_{\frp}))}
a_{\frp}(\tilde{\cax}_{\frp,z})\\
&=|\tilde{\bba}_{\frp}|\sum_{z\in\tilde{\pi}'_{\frp}(\tilde{C}''_{\frp}(\bbf_{\frp}))}
a_{\frp}(\tilde{\cax}_{\frp,z})=\frac{|\tilde{\bba}_{\frp}|}{(\caa:\cab)}
\sum_{z\in\tilde{C}''_{\frp}(\bbf_{\frp})}a_{\frp}(\tilde{\cax}''_{\frp,z''})\\
&=\frac{|\tilde{\bba}_{\frp}|}{(\caa:\cab)}q_{\frp}\frA_{\frp}(\cax'').
\end{align*}
It follows from the latter equality, Theorem \ref{avgtr1} applied to $\cax''$ and Proposition \ref{condext} applied to
$C''\to C$ that
\begin{align*}
|\frA_{\frp}(\cax')|&=\frac{|\tilde{\bba}_{\frp}|}{(\caa:\cab)}|\frA_{\frp}(\cax'')|
\le\frac{|\tilde{\bba}_{\frp}|}{(\caa:\cab)}(2g_X(2g_{C''}-2)+f_{J_{X''}})+O(q_{\frp}^{-1/2})\\
&=\frac{|\tilde{\bba}_{\frp}|}{(\caa:\cab)}\frac1{|\cab|}(2g_X(2g_{C'}-2)+f_{J_{X'}})+O(q_{\frp}^{-1/2}).
\end{align*}
\end{proof}

\section{Proof of Theorem \ref{mainthm}}

\begin{proof}
It follows from Theorem \ref{thmhp} that
$$
\rk\left(\frac{J_X(K')}{\tau B(k)}\right)=\res_{s=1}\left(\sum_{\frp\notin S}-\frA_{\frp}^*(\cax')\frac{\log
q_{\frp}}{q_{\frp}^s}\right).
$$
Note that
$$
|\frA_{\frp}^*(\cax')|\le|\frA_{\frp}(\cax')|+|a_{\frp}(B)|\le|\frA_{\frp}(\cax')|+2\dim(B)q_{\frp}^{1/2}.
$$
So, by Proposition \ref{avgtr2},
\begin{equation}\label{rkbd1}\begin{aligned}
&\rk\left(\frac{J_X(K')}{\tau B(k)}\right)\le(2g_X(2g_C-2)+f_{J_{X'}})\res_{s=1}\left(\sum_{\frp\notin
S}\frac{|\tilde{\bba}_{\frp}|}{|\caa|}\frac{\log q_{\frp}}{q_{\frp}^s}\right)\\
&+O\left(\res_{s=1}\left(\sum_{\frp\notin S}\frac{\log
q_{\frp}}{q_{\frp}^{s+1/2}}\right)\right)+2\dim(B)\res_{s=1}\left(\sum_{\frp\notin
S}\frac{\log(q_{\frp})}{q_{\frp}^{s-1/2}}\right).
\end{aligned}\end{equation}
The second series converges for $\Re(s)>1/2$, thus the corresponding residue equals 0. The second series converges for
$\Re(s)>3/2$ and can be extended meromorphically to the whole plane with just a simple pole at $s=3/2$. Hence, the
latter residue also equals 0 and so
\begin{equation}\label{rkbd1a}
\rk\left(\frac{J_X(K')}{\tau B(k)}\right)\le(2g_X(2g_C-2)+f_{J_{X'}})\res_{s=1}\left(\sum_{\frp\notin
S}\frac{|\tilde{\bba}_{\frp}|}{|\caa|}\frac{\log q_{\frp}}{q_{\frp}^s}\right).
\end{equation}

Moreover, if $\sigma\in\gal(L/k)$, then $|\tilde{\bba}_{\frp}|$ is the same for every $\frp$ such that $\sigma$ is in
the $\frp$-Frobenius conjugacy class  $(\frp,L/k)\subset\gal(L/k)$. More precisely, if $\sigma\in(\frp,L/k)$, then
$|\tilde{\bba}_{\frp}|=\#\{a\in \caa\,;\,\sigma(a)=a\}$. Denote this number by $h^0(\sigma,\caa)$. By (\ref{rkbd1a})
\begin{equation}\label{rkbd2}\begin{aligned}
&\rk\left(\frac{J_X(K')}{\tau
B(k)}\right)\le(2g_X(2g_C-2)+f_{J_{X'}})\times\\
&\left(\sum_{\tau\in\gal(L/k)}\frac{h^0(\sigma,\caa)}{|\caa|} \res_{s=1}\left(\sum_{\substack{\frp\notin S\\
\sigma\notin(\frp,L/k)}}\frac{\log q_{\frp}}{q_{\frp}^s}\right)\right).
\end{aligned}\end{equation}
By \cite{sil00} the latter residue is equal to $[L:k]=|\gal(L/k)|$. The result now follows from \cite[Lemma 9]{sil02}.
\end{proof}

\section{The rank in special towers}

In this paragraph we apply Theorem \ref{mainthm} to special towers. In the first case we assume that $C$ is an elliptic
curve defined over $k$, $C'=C$ and $\pi$ is the multiplication map $[n]:C\to C$ by an integer $n\ge1$. In particular,
$\caa$ is the subgroup $C[n]$ of $n$-torsion points of $C$. In the second case, we take $C'$ to be the pull-back of $C$
by the multiplication by $n$-map $[n]:J_C\to J_C$, where $J_C$ denotes the Jacobian variety of $C$. In particular,
$\caa$ is the subgroup $J_C[n]$ of $n$-torsion points of $J_C$.

\begin{theorem}[Serre, \cite{serre72}]\label{serrethm}
Let $C/k$ be an elliptic curve defined over a number field $k$. There is an integer $I(C/k)$ so that for every integer
$n\ge1$ the image of the representation $\rho_{C,n}:G_k\to\aut(C[n])\cong\gl_2(\bbz/n\bbz)$ has index at most $I(C/k)$
in $\aut(C[n])$.
\end{theorem}

Let $A$ be an abelian variety defined over a number field $k$ of dimension $d\ge1$. For each prime number $\ell$,
denote by $T_{\ell}(A)$ the $\ell$-adic Tate module of $A$ and by
$\rho_{\ell}:G_k\to\aut(T_{\ell}(A))\cong\gl_{2d}(\bbz_{\ell})$ the action of $G_k$ on $T_{\ell}(A)$. Let
$G_{k,\ell}=\rho_{\ell}(G_k)$ and
$\rho=\prod_{\ell}\rho_{\ell}:G_k\to\prod_{\ell}G_{k,\ell}\subset\prod_{\ell}\aut(T_{\ell}(A))$.

\begin{theorem}[Serre, \cite{serre86rib}, \cite{serre86mf}, \cite{serre86cf}]\label{serreabth}
Suppose that $k$ is large enough (depending on $A$). Then $\rho(G_k)$ is an open subgroup of $\prod_{\ell}G_{k,\ell}$.
\end{theorem}

A similar argument to \cite{serre72} shows that Theorem \ref{serreabth} is equivalent to the following result.

\begin{theorem}[Serre]\label{serreabth1}
If $k$ is large enough (depending on $A$), then exists an integer $I(A/k)\ge1$ such that for every integer $n\ge1$ the
image of the representation $\rho_{A,n}:G_k\to\aut(A[n])\cong\gl_{2d}(\bbz_{\ell})$ has index at most $I(A/k)$ in
$\aut(A[n])$.
\end{theorem}

\begin{theorem}\label{thmabtow}
Assume one of the two following situations holds:
\begin{enumerate}
\item[(a)] $C$ is an elliptic curve, $C_n=C'=C$, $\pi$ is equal to $[n]:C\to C$.

\item[(b)] $C_n=C'$ is the pull-back of $C$ by $[n]:J_C\to J_C$.
\end{enumerate}
Let $\cax_n$ be the pull-back of $\cax$ by $\pi$, $K_n=k(C_n)$ and $X_n$ the generic fiber of $\phi_n:\cax_n\to C_n$.
Assume that Tate's conjecture is true for $\cax_n/k$. Then (for each of the above cases):
\begin{enumerate}
\item For every integer $n\ge1$:
\begin{enumerate}
\item[(i-a)]
$$
\rk\left(\frac{J_X(K_n)}{\tau B(k)}\right)\le f_{J_{X_n}}I(C/k)\frac{d(n)}{n^2},
$$
where $d(n)$ denotes the number of positive divisors of $n$.

\item[(i-b)] If $k$ is large enough (depending on $J_X$), denoting by $I(J_C/k)$ the constant of Theorem
\ref{serreabth1} corresponding to the Jacobian variety $J_C$ of $C$, then
$$
\rk\left(\frac{J_X(K_n)}{\tau B(k)}\right)\le(2g_X(2g_{C_n}-2)+f_{J_{X_n}})I(J_C/k)\frac{d(n)}{n^{2g_C}}.
$$
\end{enumerate}
\item In both cases, the sum
$$
\frac1x\sum_{n\le x}\frac1{\log(f_{J_{X_n}})}\rk\left(\frac{J_X(K_n)}{\tau B(k)}\right)
$$
is bounded as $x\to\infty$. Thus the average rank of $J_X(K_n)/\tau B(k)$ is smaller than a fixed multiple of the
logarithmic of its conductor.

\item In both cases, there exists a constant $\kappa=\kappa(k,C,J_X)$ (resp. $\kappa=\kappa(k,J_C,$ \linebreak $J_X)$)
so that for sufficiently large $n$ we have
$$
\rk\left(\frac{J_X(K_n)}{\tau B(k)}\right)\le f_{J_{X_n}}^{\kappa/\log(\log(f_{J_{X_n}}))}.
$$
In particular, for every $\epsilon>0$ we have
$$
\rk\left(\frac{J_X(K_n)}{\tau B(k)}\right)\ll f_{J_{X_n}}^{\epsilon},
$$
where the implied constant depends on $k$, $C$ (resp. $J_C$), $J_X$ and $\epsilon$, but not on $n$.
\end{enumerate}
\end{theorem}

\begin{proof}
In the first case, by Theorem \ref{serrethm} and \cite[Lemma 10]{sil02},
$$
\#\frO_{G_k}\le I(C/k)\#\frO_{\aut(C[n])}(C[n])=I(C/k)\#\frO_{\gl_2(\bbz/n\bbz)}((\bbz/n\bbz)^2).
$$
It follows from \cite[Proposition 15]{sil02} that
$$
\#\frO_{\gl_2(\bbz/n\bbz)}((\bbz/n\bbz)^2)=d(n).
$$
Hence (i-a) follows from Theorem \ref{mainthm}. In the second case, we apply Theorem \ref{serreabth1} and \cite[Lemma
10]{sil02} to obtain
\begin{align*}
\#\frO_{G_k}(J_C[n])&\le I(J_C/k)\#\frO_{\aut(J_C[n])}(J_C[n])\\
&=I(J_C/k)\#\frO_{\gl_{2g_C}(\bbz/n\bbz)}((\bbz/n\bbz)^{2g_C}).
\end{align*}
Again \cite[Proposition 15]{sil02} implies
$$
\#\frO_{\gl_{2g_C}(\bbz/n\bbz)}((\bbz/n\bbz)^{2g_C})=d(n).
$$
So (i-b) follows from Theorem \ref{mainthm}.

By Proposition \ref{condext}, in the first case (resp. in the second case), $f_{J_{X_n}}=n^2f_{J_X}$ (resp.
$f_{J_{X_n}}=n^{2g_C}f_{J_X}$) and $2g_{C_n}-2=n^{2g_C}(2g_C-2)$ in the second case. It follows from \cite[Theorem
3.3]{apo76} that the function $d(n)$ satisfies the following property
$$
\sum_{n\le x}d(n)\sim x\log(x).
$$
Thus,
$$
\frac1x\sum_{2\le n\le x}\frac{d(n)}{\log(n)}
$$
is bounded for all $x\ge2$. Therefore,
\begin{align*}
&\frac1x\sum_{n\le x}\frac1{\log(f_{J_{X_n}})}\rk\left(\frac{J_X(K_n)}{\tau B(k)}\right)\\
&\le\frac1x\sum_{n\le x}\frac1{\log(n^2f_{J_X})}(f_{J_X}d(n)I(C/k))
\end{align*}
is also bounded (in the first case), as well as
\begin{align*}
&\frac1x\sum_{n\le x}\frac1{\log(f_{J_{X_n}})}\rk\left(\frac{J_X(K_n)}{\tau B(k)}\right)\\
&\le\frac1x\sum_{n\le x}\frac1{\log(n^{2g_C}f_{J_X})}((2g_X(2g_C-2)+f_{J_X})d(n)I(C/k))
\end{align*}
(in the second case). Whence (2) follows.

Finally, item (3) follows as in \cite[Theorem 16]{sil02}.
\end{proof}

\begin{remark}
In the two types of unramified abelian towers considered in this paragraph, Theorem \ref{thmabtow} (3) shows that the
rank grows slower along the tower than the geometric bound.

In \cite[Th\'eor\`eme 2.1]{mic97} Michel computed (under the assumption of the validity of the standard conjectures) an
upper bound for the average rank of a family of abelian varieties over $\bbq$. This bound depended (among other
invariants) on the average logarithmic conductor. This type of problem is ``horizontal'' (fix the field and vary the
abelian variety), whereas ours (eg. Theorem \ref{thmabtow} (2)) is ``vertical'' (fix the abelian variety and vary the
field).
\end{remark}

\end{document}